\documentclass[12pt]{amsart}
\usepackage{amsfonts, amssymb, amsmath, amsthm, latexsym, array}
\usepackage{fullpage}
\usepackage{verbatim}

\newtheorem{thm}{Theorem}
\newtheorem{lm}{Lemma}

\theoremstyle{remark}
\newtheorem{ex}{Example}

\theoremstyle{definition}

\newcommand{\cle }{\noindent {\bf Corollary.} }

\def\Hom{\mathrm{Hom}}
\def\ad{\mathrm{ad}}
\def\diag{\mathrm{diag}}
\def\Ad{\mathrm{Ad}}

\def\SO{{\rm SO}}

\newcommand{\rk}{{\rm rk}}
\newcommand{\gt}{\mathfrak}
\newcommand{\id}{{\rm id}}
\newcommand{\GL}{{\rm GL}}
\newcommand{\Lie}{{\rm Lie\,}}

\newcommand {\Ker}{{\mathrm{Ker\,}}}
\newcommand {\g}{\mathfrak{g}}
\newcommand {\fN}{\mathfrak{N}}

\newfam\Bbbfam\newfam\eufam\newfam\eusfam%
\font\Bbbfont=msbm10 scaled 1200%
\font\Bbbsmallfont=msbm8%
\textfont\Bbbfam=\Bbbfont\scriptfont\Bbbfam=\Bbbsmallfont
\font\euzw=eufm10 scaled 1200%
\font\euac=eufm7 scaled 1200%
\font\euacc=eufm7 scaled 1000%
\textfont\eufam=\euzw\scriptfont\eufam=\euac%
\scriptscriptfont\eufam=\euacc%
\font\euszw=eusm10 scaled 1200%
\font\eusac=eusm7 scaled 1200%
\font\eusacc=eusm7 scaled 1000%
\textfont\eusfam=\euszw\scriptfont\eusfam=\eusac%
\scriptscriptfont\eusfam=\eusacc%

\def\bbk{\hbox {\Bbbfont\char'174}}

\begin{document}
\title[(Ir)Reducibility of commuting varieties]{(Ir)Reducibility of
some commuting \\ varieties  associated with involutions}
\author[O.\,Yakimova]{\sc Oksana Yakimova}
\maketitle
\section*{Introduction}

The ground field $\bbk$ is algebraically closed and of characteristic zero.
Let $\gt g$ be a reductive algebraic  Lie algebra over $\bbk$
and $\sigma$ an involutory automorphism of $\gt g$.
Then $\gt g=\gt g_0\oplus\gt g_1$ is the direct sum of $\sigma$-eigenspaces.
Here $\gt g_0$ is a reductive subalgebra and
$\gt g_1$ is a $\gt g_0$-module.
Let $G$ be the adjoint group of $\gt g$
and $G_0\subset G$ a connected subgroup with $\Lie G_0=\gt g_0$.
The {\it commuting variety} of
$(\mathfrak{g},\gt g_0)$ is the following set:
$$
{\gt C}(\gt g_1) = \{ (x,y) \in \gt g_1 \times \gt g_1 \  \mid \ [x,y] = 0
\}.
$$
The problem whether ${\gt C}(\gt g_1)$
is irreducible was considered by Panyushev
\cite{Dima1}, \cite{Dima2} and Sabourin-Yu \cite{sy}, \cite{sy2}.
Suppose $\gt g$ is simple. Then the known results are
\begin{itemize}
\item[$\bullet$] if the rank of the symmetric pair
$(\mathfrak{g},\gt g_0)$ is equal to the
semisimple rank of $\mathfrak{g}$ (called the maximal rank case), then
the corresponding commuting variety is irreducible, \cite{Dima1};
\item[$\bullet$]  if the rank of $(\mathfrak{g},\gt g_0)$ equals
$1$, then ${\gt C}(\gt g_1)$ is irreducible only in one case, namely,
$(\gt{so}_{m+1}, \gt{so}_{m})$, \cite{Dima2}, \cite{sy};
\item[$\bullet$] for
$(\gt{sl}_{2n}, \gt{sp}_{2n})$ and $(E_6,F_4)$ the
corresponding commuting variety is irreducible, \cite{Dima2};
\item[$\bullet$] if $(\gt g,\gt g_0)=(\gt{so}_{2+m},
\gt{so}_2{\oplus}\gt{so}_m)$, then ${\gt C}(\gt g_1)$
is irreducible, \cite{sy2}.
\end{itemize}
For all other symmetric pairs the problem is open. In
sections~1--3, we extend the result of \cite{sy2} to all
symmetric pairs $(\gt{so}_{n+m}, \gt{so}_n\oplus\gt{so}_m)$.
The scheme of the proof is similar to that of \cite{sy2}.
But as it often happens, the argument in a general situation
is shorter and simpler, than in a particular case.

In \cite{Dima2}, it was conjectured that ${\gt C}(\gt g_1)$
is irreducible if the rank of the symmetric pair is greater
than $1$. This conjecture is not true. In section~4, we
prove that $\gt C(\gt g_1)$ is reducible for symmetric pairs
$(\gt{gl}_{n+m},\gt{gl}_n{\oplus}\gt{gl}_m)$ with $n\ne m$,
$(\gt{so}_{2n},\gt{gl}_n)$ with odd $n$, and $(E_6,
\gt{so}_{10}{\oplus}\bbk)$.

\section{Commuting variety}

The commuting variety ${\gt C}(\gt g)=\{(x,y)\in\gt g+\gt
g\mid[x,y]=0\}$ of a reductive Lie algebra $\gt g$ was
considered by  Richardson in \cite{r}, where he shows that
${\gt C}(\gt g)$ coincides with the closure of ${G(\gt
t{\times}\gt t)}$ for any  Cartan subalgebra $\gt
t\subset\gt g$, and is therefore irreducible. It is not yet
known whether ${\gt C}(\gt g)$ is normal or whether the
ideal of ${\gt C}(\gt g)$ is generated by quadrics.

Let $(\gt g,\gt g_0)$ be a symmetric pair and
$\mathfrak{c}\subset\gt g_1$ a maximal abelian subspace
consisting of semisimple elements. Any such subspace is
called a {\it Cartan subspace} of $\gt g_1$. All Cartan
subspaces are $G_0$-conjugate, see \cite{kr}. The dimension
of $\gt c$ is called the {\it rank} of the symmetric pair
$(\gt{g}, \gt g_0)$. It is well-known that ${\gt C}_{0} =
\overline{G_0(\gt{c}{\times}\gt{c})}$ is the unique
irreducible component of ${\gt C}(\gt g_1)$ of maximal
dimension, see \cite[Sect. 3]{Dima1}. Here $\dim{\gt C}(\gt
g_1)=\dim\gt g_1 + \dim\gt{c}$. It follows that $\gt C(\gt
g_1)$ is irreducible if and only if $\gt C(\gt g_1)=\gt
C_0$.

The irreducibility problem for
$\gt C(\gt g_1)$ was first considered by Panyushev in
\cite{Dima1}. As was mentioned above, $\gt C(\gt g_1)$ is not always
irreducible. On the other hand, in some particular cases
one can say more about properties of $\gt C(\gt g_1)$.
If $(\gt g,\gt g_0)$ is a symmetric pair of
maximal rank, then $\gt C(\gt g_1)$ is a normal
complete intersection and the ideal of $\gt C(\gt g_1)$
in $\bbk[\gt g_1\times\gt g_1]$
is generated by quadrics, see \cite{Dima1}.

Let $n$ be a non-negative integer. The set
$$
\gt g_1^{(n)}=\{\xi\in \gt g_1\mid \dim G_0\xi=n \}
$$
is locally closed. Irreducible components
of $\gt g_1^{(n)}$ are called {\it $G_0$-sheets} of $\gt g_1$.
The following statement is well-known for the specialists, but
we include a proof here for the sake of completeness.

\begin{lm}\label{sheet}
Let $S$ be a $G_0$-sheet of $\gt g_1$ containing semisimple
elements. Suppose for each semisimple $h\in S$ we have
$(\{h\}{\times}(\gt g_1)_h)\subset\gt C_0$. Then
$(\{x\}{\times}(\gt g_1)_x)\subset\gt C_0$ for each $x\in
S$.
\end{lm}
\begin{proof}
Let $x\in S$. Since $S$ contains semisimple
elements, they form a dense subset. Therefore,
we can find a morphism
$\gamma:\bbk\to S$ such that
$\gamma(0)=x$ and $\gamma(t)$ is semisimple
for each $t\ne 0$. Then
$(\gt g_1)_x=\lim_{t\to 0}(\gt g_1)_{\gamma(t)}$, where the
limit is taken in an appropriate Grassmannian.
For each $y\in(\gt g_1)_x$, we can define
elements $y(t)\in(\gt g_1)_{\gamma(t)}$ such that
$y=\lim_{t\to 0}y(t)$.
Since $(x,y)=\lim_{t\to 0}(\gamma(t),y(t))$ and
$(\gamma(t),y(t))\in\gt C_0$ for each $t\ne 0$,
we conclude that $(x,y)\in\gt C_0$.
\end{proof}

\section{Semisimple and nilpotent elements in $\gt g_1$}
\label{centralisers}

In this section $\gt g=\gt{so}_{n+m}$, $\gt g_0=\gt{so}_n{\oplus}\gt{so}_m$.
Let $V=\bbk^{n+m}$ be
a vector space of
the defining representation of $\gt g$.
Then we have a $G_0$-invariant decomposition
$V=V_a{\oplus}V_b$, where
$V_a=\bbk^n$, $V_b=\bbk^m$, and
$\gt g_1\cong\bbk^n{\otimes}\bbk^m$ as a $G_0$-module.
Denote by $(\,\,,\,)$ the non-degenerate symmetric
$\gt g$-invariant bilinear form on $V$.

Let $\gt g_h$ be the centraliser of an element $h\in\gt g_1$.
Then $\sigma$ induces the symmetric decomposition
$(\gt g_h)=(\gt g_0)_h\oplus(\gt g_1)_h$,
where $(\gt g_0)_h$ is the centraliser of
$h$ in $\gt g_0$. First we describe centralisers
$\gt g_h$ and $(\gt g_0)_h$ of semisimple elements $h\in\gt g_1$.

\begin{lm}\label{ss}
Let $h\in\gt g_1$ be a semisimple element.
Then the symmetric pair $(\gt g_h, (\gt g_0)_h)$ is a direct sum
$(\bigoplus_{i=1}^{r}(\gt{gl}_{k_i},\gt{so}_{k_i}))\oplus
(\gt{so}_{n+m-2k},\gt{so}_{n-k}{\oplus}\gt{so}_{m-k})$,
where $k=\sum_i k_i$.
\end{lm}
\begin{proof}
Recall several well-known facts about semisimple elements
of $\gt{so}_{n+m}$. Let
$v_\lambda\in V$ be an eigenvector of $h$ such that
$h\cdot v_\lambda=\lambda v_\lambda$ and $\lambda\ne 0$.
Since $h$ preserves the symmetric form
$(\,\,,\,)$, we have $(v_\lambda, v_\lambda)=0$.
Also if $h\cdot v=\lambda v$, $h\cdot w=\mu w$, then
$(v,w)\ne 0$ only if $\lambda=-\mu$.
Let $\{\pm\lambda_i, 0\mid i=1,\ldots, r\}$
be the set of the eigenvalues of $h$.
Then there is an orthogonal
$h$-invariant decomposition
$$
V=(V_{\lambda_1}{\oplus}V_{-\lambda_1})\oplus\ldots\oplus
(V_{\lambda_r}{\oplus}V_{-\lambda_r})\oplus V_0.
$$
Here each $V_{\lambda_i}$ is an isotropic subspace,
$(V_{\lambda_i}, V_{\pm\lambda_j})=0$ if
$\lambda_i\ne\pm\lambda_j$ and $(V_0,V_{\pm\lambda_i})=0$
for each $\lambda_i$. Therefore
$\gt g_h\subset(\bigoplus_{i=1}^r \gt{so}(V_{\lambda_i}{\oplus}V_{-\lambda_i}))
\oplus\gt{so}(V_0)$. More precisely,
if $\dim V_{\lambda_i}=k_i$ and $k=\sum_i k_i$, then
$\gt g_h=(\bigoplus_i\gt{gl}_{k_i})\oplus\gt{so}_{n+m-2k}$.

Now it remains to describe $(\gt g_0)_h=(\gt g_h)^\sigma$.
We may assume that $\sigma$ is a conjugation by
a diagonal matrix $A\in{\rm O}_{n+m}$
such that $A|_{V_a}=-\id$ and $A|_{V_b}=\id$.
Since $\sigma(h)=-h$, we have
$A\cdot V_{\lambda_i}=V_{-\lambda_i}$ and $A\cdot V_0=V_0$.
Moreover,
$A$ determines a non-degenerate symmetric form $(\,\,,\,)_A$
on each $V_{\lambda_i}$ by the formula
$(v,w)_A=(v,A\cdot w)$.
Therefore, each $\gt{so}(V_{\lambda_i}{\oplus}V_{-\lambda_i})$
is $\sigma$-invariant,
$(\gt{so}(V_{\lambda_i}{\oplus}V_{-\lambda_i}))^\sigma\cong
\gt{so}_{k_i}{\oplus}\gt{so}_{k_i}$, and
$(\gt{gl}_{k_i})^\sigma=\gt{so}_{k_i}$. Finally,
the restriction $A|_{V_0}$ has signature $(n-k,m-k)$.
Thus $(\gt{so}_{n+m-2k})^\sigma=\gt{so}_{n-k}{\oplus}\gt{so}_{m-k}$.
\end{proof}

Denote by $\fN(\gt g_1)$ the {\it nullcone\/} of $G_0:\gt
g_1$, i.e., the set of all nilpotent elements in $\gt g_1$.
Recall several standard facts concerning nilpotent elements
in $\gt{gl}(V)$. Suppose $e\in\gt{gl}(V)$ is nilpotent and
$m=\dim \Ker(e)$. Then by the theory of Jordan normal form,
there are vectors $w_1,\ldots,w_m\in V$ and non-negative
integers $d_1,\ldots,d_m$ such that $e^{d_i+1}{\cdot}w_i=0$
and $\{ e^s{\cdot}w_i \mid 1\le i\le m, \ 0\le s\le d_i\}$
is a basis for $V$. Let $V_i\subset V$ be a linear span of
$\{w_i,e{\cdot}w_i,\ldots, e^{d_i}{\cdot}w_i\}$. Then the
spaces $\{V_i\}$ are called the Jordan (or {\it cyclic})
spaces of the nilpotent element $e$ and $V=\oplus_{i=1}^m
V_i$.

\begin{lm}\label{ab}
Suppose $e\in\fN(\g_1)$.  Then
the cyclic vectors $\{w_i\}_{i=1}^m$ and hence the cyclic spaces
$\{V_i\}$'s can be chosen such that the following properties are satisfied:
\begin{itemize}
\item[\sf (i)] \
there is an involution $i\mapsto i^*$ on the set
$\{1,\dots, m\}$ such that $d_i=d_{i^*}$,
$i=i^*$ if and only if $\dim V_i$ is odd, and
$(V_i, V_j)=0$ if $i\ne j^*$;
\item[\sf (ii)] \ $\sigma(w_i)=\pm w_i$.
\end{itemize}
\end{lm}
\begin{proof}
Part (i) is a standard property of the nilpotent orbits in
$\gt{so}(V)$, see, for example, \cite[Sect.~5.1]{cm} or \cite[Sect.~1]{ja}.
Then part (ii) says that in the presence of the involution $\sigma$
cyclic vectors for $e\in\fN(\g_1)$ can be chosen to be $\sigma$-eigenvectors,
see \cite[Prop.~2]{ohta}.
\end{proof}

For each $e\in\fN(\gt g_1)$ we choose cyclic vectors
$\{w_i\}$ as prescribed by Lemma~\ref{ab}. Say that
$e^s{\cdot}w_i$ has type $a$ if
$\sigma(e^s{\cdot}w_i)=-e^s{\cdot}w_i$, i.e.,
$e^s{\cdot}w_i\in V_a$; and  $e^s{\cdot}w_i$ has type $b$ if
$e^s{\cdot}w_i\in V_b$. Since $\sigma(e)=-e$, if
$e^s{\cdot}w_i\in V_a$, then $e^{s+1}{\cdot}w_i\in V_b$ and
vice versa. Therefore each string $\langle
w_i,e{\cdot}w_i,\ldots,e^{d_i}{\cdot}w_i\rangle$
has one of the following types: \\
\phantom{$n$} \qquad\quad
$aba\ldots ab$, \qquad  $bab\ldots ba$, \qquad $aba\ldots ba$,
\qquad $bab\ldots ab$.

Let $e\in\fN(\gt g_1)$. There is an
$\gt{sl}_2$-triple $\{e,f,h\}$ such that
$f\in\gt g_1$ and $h\in\gt g_0$.
Recall that $e$ is called {\it even} if the eigenvalues of
$\ad(h)$ on $\gt g$ are even. An element
$e\in\fN(\gt g_1)$ is said to be
{\it $\sigma$-distinguished} (in other notations
$\gt p$- or $({-}1)$-distinguished) if $(\gt g_1)_e$ contains
no semisimple elements of $[\gt g,\gt g]$.

\begin{lm}\label{ev}
In case $(\gt g,\gt g_0)=(\gt{so}_{n+m}, \gt{so}_n{\oplus}\gt{so}_m)$
each $\sigma$-distinguished element $e\in\fN(\gt g_1)$
is even.
\end{lm}
\begin{proof}
Let $\{V_i\}$ be the cyclic spaces of $e\in\fN(\gt g_1)$
chosen as prescribed by Lemma~\ref{ab}.
Suppose there is an even-dimensional $V_i$.
According to \cite[Prop.~2]{ohta}, if
$V_i$ has type $aba\ldots ab$, then
$V_{i^*}$ has type $bab\ldots ba$,
i.e., if $\sigma(w_i)=-w_i$,
then $\sigma(w_{i^*})=w_{i^*}$.
Let $\gt l$ be a Levi subalgebra of $(\gt{so}(V_i{\oplus}V_{i^*}))_e$.
We may assume that $\gt l$ is $\sigma$-invariant.
Then $\gt l=\gt{sl}(\bbk w_i{\oplus}\bbk w_{i^*})\cong\gt{sl}_2$.
The restriction of $\sigma$ defines a symmetric
decomposition $\gt l=\gt l_0\oplus\gt l_1$, where
$\gt l_0=\gt l^\sigma=\gt{so}_2$. Therefore
$\gt l_1=\gt l\cap\gt g_1$ contains semisimple
elements. This means that
$e$ is not $\sigma$-distinguished. Hence, all
$V_i$ are odd-dimensional and $e$ is even.
\end{proof}

\section{${\gt C}(\gt g_1)$ is irreducible in case of
$(\gt{so}_{n+m},\gt{so}_n{\oplus}\gt{so}_m)$}\label{f}

In this section we prove that
$\gt C(\gt g_1)$ is irreducible in case
$(\gt g,\gt g_0)=(\gt{so}_{n+m},\gt{so}_n{\oplus}\gt{so}_m)$.

The following lemma is taken from \cite{sy2}, but the proof
given here is shorter. Note that this lemma is valid for any symmetric pair
$(\gt g,\gt g_0)$.

\begin{lm}\label{even}
Suppose $e\in\fN(\gt g_1)$ is even. Then
$e$ belongs to a $G_0$-sheet containing semisimple elements.
\end{lm}
\begin{proof}
Let $(e,f,h)$ be an $\gt{sl}_{2}$-triple
such that $f\in\gt g_1$, $h\in\gt g_0$.
Since $e$ is even, we have $\dim\gt g_h=\dim\gt g_e$.
Set $e(t):=e-t^2f$  for $t\in\bbk$.
If $t\ne 0$, then $e(t)$ is semisimple and conjugated
to $th$. Therefore $\dim\gt g_{e(t)}=\dim\gt g_h=\dim\gt g_e$.
Clearly $e(0)=e=\lim_{t\to 0} e(t)$ and
the $G_0$-sheet containing $e$ contains also semisimple elements $e(t)$.
\end{proof}

Suppose $(\gt g,\gt g_0)=(\gt{so}_{n+m},\gt{so}_n{\oplus}\gt{so}_m)$ and
let $\gt{c}\subset\gt g_1$ be a Cartan subspace.

\begin{thm}\label{reduction}
The commuting variety ${\gt C}(\gt g_1)$ is irreducible.
\end{thm}
\begin{proof}
Recall that ${\gt C}_{0} =\overline{G_0(\gt c{\times}\gt
c)}$. Following the original proof of Richardson \cite{r}
(see also \cite[Sect.~2]{Dima2}), we show by induction on
$\dim\gt c$ that $\gt C(\gt g_1)=\gt C_0$. The base of
induction is the rank~$1$ case $(\gt{so}_{n+1},\gt{so}_n)$,
where the irreducibility of $\gt C(\gt g_1)$ is proved in
\cite{Dima2}, \cite{sy}. Let $(x,y) \in {\gt C}(\gt g_1)$.

(1) Suppose
there is a  semisimple element $h\in\gt g_1$
such that $[h,x]=[h,y]=0$.
This assumption is automatically satisfied if either $x$ or $y$
is semisimple. Moreover, if $x$ (or $y$) is not nilpotent and
$x=x_s+x_n$ is the Jordan decomposition, then
$x_s\in\gt g_1$ and $[x_s,x]=0$, $[x_s,y]=0$.

Consider the symmetric pair $(\gt g_h, (\gt g_0)_h)$.
Replacing $\gt c$ by a conjugated Cartan subspace, we may
assume that $h\in\gt c$. Then $\gt c$ is a Cartan subspace
of $(\gt g_1)_h$. Also, $x,y\in(\gt g_1)_h$ by the
assumption. By Lemma~\ref{ss}, $(\gt g_h, (\gt g_0)_h)=
(\bigoplus_{i=1}^{r}(\gt{gl}_{k_i},\gt{so}_{k_i}))\oplus
(\gt{so}_{n+m-2k},\gt{so}_{n-k}{\oplus}\gt{so}_{m-k})$. Note
that each $(\gt{gl}_{k_i},\gt{so}_{k_i}))$ is a symmetric
pair of maximal rank, hence, the corresponding commuting
variety is irreducible, see \cite[(3.5)(1)]{Dima1}. Clearly,
the commuting variety corresponding to a direct sum of
symmetric pairs is a direct product of the commuting
varieties corresponding to the summands. Therefore, using
the inductive hypothesis, we conclude that $\gt C((\gt
g_1)_h)$ is irreducible. Thus $(x,y)\in\overline{(G_0)_h(\gt
c{\times}\gt c)}$ and, hence, $(x,y)\in\gt C_0$.

(2) It remains to consider pairs of commuting nilpotent elements.
Suppose first that there is a semisimple element
$h\in\gt g_1$ such that $[x,h]=0$.
Then $(x,(1-t)y+th)\in\gt C(\gt g_1)$ for each
$t\in\bbk$ and $(1-t)y+th$ is nilpotent only for
a finite number of $t$'s. Therefore, by part~(1),
one has $(x,(1-t)y+th)\in\gt C_0$ for almost all $t$.
Since $y=\lim_{t\to 0}(1-t)y+th$, we get $(x,y)\in\gt C_0$.

(3) Now we may assume that both $x$ and $y$ are
$\sigma$-distinguished nilpotent elements.
According to Lemma~\ref{ev}, $x$ is even. Then,
by Lemma~\ref{even}, $x$ belongs to a
$G_0$-sheet containing semisimple elements.
According to part~(1), the assumptions of
Lemma~\ref{sheet} are satisfied and it follows that
$(x,y)\in\gt C_0$.
\end{proof}

\section{Several new examples of reducible commuting varieties}

Let $\gt g=\gt g(-1)\oplus\gt g(0)\oplus\gt g(-1)$ be a
short grading of a Lie algebra $\gt g$. Set $\gt g_0:=\gt
g(0)$. Then $(\gt g,\gt g_0)$ is a symmetric pair with
$\gt g_1=\gt g(-1){\oplus}\gt g(1)$. 
Let $\gt c\subset\gt g_1$ be a Cartan subspace and $\gt c(-1)$, $\gt
c(1)$ the images of $\gt c$ under projections on $\gt g(-1)$
and $\gt g(1)$, respectively.

\begin{lm}\label{piki}
Suppose $\gt C(\gt g_1)$ is irreducible. Then
\begin{center}
$({\bf \spadesuit})\quad$ $\overline{G_0(\gt c(1){\times}\gt
c(1))}=\gt g(1){\times}\gt g(1)$.
\end{center}
\end{lm}
\begin{proof}
Since each $\gt g(-1)$ and $\gt g(1)$ consists of nilpotent
elements, we have $\dim\gt c({\pm 1})=\dim\gt c$ and $\gt
c\cong\diag(\gt c({-}1){\oplus}\gt c(1))$. Clearly $\gt
g(1){\times}\gt g(1)\subset\gt C(\gt g_1)$. Since $\gt C(\gt
g_1)=\gt C_0=\overline{G_0(\gt c{\times}\gt c)}$, we get
$\gt g(1){\times}\gt g(1)=\gt g(1){\times}\gt g(1)\cap\gt
C_0=\overline{G_0(\gt c(1){\times}\gt c(1))}$.
\end{proof}

\vskip0.2ex
\cle {\it If condition $({\bf \spadesuit})$ is not
satisfied, then $\gt C(\gt g_1)$ has at least three
irreducible components.}

\vskip0.2ex

Now we give three examples of symmetric pairs arising from
short gradings such that condition $({\bf \spadesuit})$ is
not satisfied for them.

\begin{ex} Suppose $(\gt g, \gt
  g_0)=(\gt{gl}_{n+m},\gt{gl}_n\oplus\gt{gl}_m)$.
Let $V$ be a $(n+m)$-dimensional vector space such that $\gt
g=\gt{gl}(V)$. Let $V=V_a\oplus V_b$ be the $\gt
g_0$-invariant decomposition with $\dim V_a=n$, $\dim
V_b=m$. Then the involution $\sigma$ is induced by a
diagonal matrix $A\in\GL(V)$ such that $A|_{V_a}=-\id$,
$A|_{V_b}=\id$. We have $\gt g(1)= \Hom(V_a,V_b)$, $\gt
g(-1)=\Hom(V_b,V_a)$ and $\gt g_1=\gt g(1)\oplus\gt g(-1)$.
Assume that $n\le m$. Suppose $\xi,\eta\in\gt g_1$ and
$[\xi,\eta]=0$. If $\xi=X+Y$, $\eta=Z+U$, where
$X,Z\in\Hom(V_a,V_b)$, $Y,U\in\Hom(V_b,V_a)$, then we set
$D_1(\xi,\eta):=(X|Z)$, 
where $(X|Z)$ is an $m{\times}2n$ matrix.

Now suppose $n\ne m$, i.e., $n<m$. We show that in this case
condition $({\bf \spadesuit})$ is not satisfied.
Let $\gt c=\{X+X^t\mid X=(x_{i,j}), x_{i,j}=0 \mbox{ if } i\ne j\}\subset\gt g_1$
be a Cartan subspace. Then $\gt c(1)=\{X=(x_{i,j}) \mid x_{i,j}=0 \mbox{ if } i\ne j\}$.
It is clear, that $\rk D_1(t,h)\le n$ for each pair
$(t,h)\in\gt c{\times}\gt c$. Let $g\in G_0$. Then
$g=B\times C$, where $B\in\GL_n=\GL(V_a)$,
$C\in\GL_m=\GL(V_b)$. If $\xi=X+Y\in\gt g_1$, then
$\Ad(g){\cdot}\xi=CXB^{-1}+BYC^{-1}$. Therefore,
$g{\cdot}D_1(\xi,\eta):=D_1(\Ad(g){\cdot}\xi,\Ad(g){\cdot}\eta)=
CD_1(\xi,\eta)\hat B$, where
$$
\hat B=\left(\begin{array}{cc}
B^{-1} & 0\\
0 & B^{-1} \\
\end{array}
\right)
$$
is a non-degenerate $2n{\times}2n$ matrix.

Since $\rk
D_1(\Ad(g){\cdot}\xi,\Ad(g){\cdot}\eta)=\rk D_1(\xi,\eta)$,
we have 
$\rk D_1(t,h)\le n$ for each pair $(t,h)\in\gt
C_0=\overline{G_0(\gt c{\times}\gt c)}$.

It remains to find a pair of  matrices $(X,Y)\in\gt
g(1){\times}\gt g(1)$ such that $\rk (X|Y)>n$. Set
$X:=(x_{ij})$ where $x_{ij}=1$ if $i=j$, and $x_{ij}=0$
otherwise; $Y:=(y_{ij})$ where $y_{ij}=1$ if $i=j+m-n$ and
$y_{ij}=0$ otherwise. It is easy to see, that $(X|Y)$ is a
matrix of the maximal possible rank, which equals
$\min(m,2n)$. In particular, $\rk (X|Y)>n$ and
$(X,Y)\not\in\overline{G_0(\gt c(1){\times}\gt c(1))}$.
\end{ex}

\begin{ex} Let $\gt g=\gt{so}(V)=\gt{so}_{2n}$,
where $2n=\dim V$. Consider a decomposition of $V$ into a
direct sum of two isotropic subspaces $V=V_+\oplus V_-$.
Suppose $A\in{\rm O}(V)$ and $A|_{V_+}=\id$,
$A|_{V_-}=-\id$. Then a conjugation by $A$ defines an
involution of $\gt g$ such that $\gt
g_0=\gt{gl}_n\cong\gt{gl}(V_+)\cong\gt{gl}(V_-)$. We have
$\gt g_1\subset\Hom(V_+,V_-)\oplus\Hom(V_-,V_+)$. More
presicely, one can choose a basis in $V_+\oplus V_-$ such
that
$$
G_0=\left\{\left(\begin{array}{cc}
          B & 0 \\
          0 & (B^t)^{-1} \\
         \end{array}\right)\mid B\in\GL_n \right\}, \enskip
\gt g_1=\left\{\left(\begin{array}{cc}
         0 & X \\
         Y & 0 \\
         \end{array}\right) \mid X,Y\in \gt{gl}_n,
 X=-X^t, Y=-Y^t \right\}.
$$
Here $\gt g(1)=\gt g_1\cap\Hom(V_+,V_-)\cong\gt{so}_n$, $\gt
g(-1)=\gt g_1\cap\Hom(V_-,V_+)\cong\gt{so}_n$. For each pair
$(\xi,\eta)\in\gt g_1\times\gt g_1$ with
$\xi=\left(\begin{array}{cc}
         0 & X \\
         Y & 0 \\
         \end{array}\right)$, $\eta=
\left(\begin{array}{cc}
         0 & Z \\
         U & 0 \\
         \end{array}\right)$
we set $D_1(\xi,\eta):=(X|Z)$. 
Take a Cartan subspace $\gt c\subset\gt g_1$ consisting of
skew symmetric anti-diagonal matrices, i.e.,
$$\gt c
=\left\{\left(
\begin{array}{cc}
         0 & X \\
         X & 0 \\
\end{array}\right) \mid X=(x_{i,j}), \, x_{ij}=0 \mbox{ if }
i\ne(n+1-j) \right\}.$$

Now suppose $n$ is odd and $n=2k+1$. When for each pair
$(t,h)\in\gt c{\times}\gt c$ we have $\rk D_1(t,h)<n$.
Suppose $g=\left(\begin{array}{cc}
          B & 0 \\
          0 & (B^t)^{-1} \\
         \end{array}\right)\in G_0$.
If $\xi=\left(\begin{array}{cc}
         0 & X \\
         Y & 0 \\
         \end{array}\right)\in\gt g_1$, then
$\Ad(g){\cdot}\xi=\left(\begin{array}{cc}
         0 & BXB^t \\
         (B^t)^{-1}YB^{-1} & 0 \\
         \end{array}\right)$. Therefore,
$g{\cdot}D_1(\xi,\eta):=D_1(\Ad(g){\cdot}\xi,\Ad(g){\cdot}\eta)=
BD_1(\xi,\eta)\hat B$, where $ \hat
B=\left(\begin{array}{cc}
B^t & 0\\
0 & B^{t} \\
\end{array}
\right)$ is a non-degenerate $2n{\times}2n$ matrix.  Hence,
$\rk D_1(\Ad(g){\cdot}\xi,\Ad(g){\cdot}\eta)=\rk
D_1(\xi,\eta)$ and $\rk D_1(t,h)< n$ for each pair
$(t,h)\in\gt C_0$.

Let $X,Z\in\gt g(1)$ be skew-symmetric $n{\times}n$ matrices
of rank $2k$ such that the last column and the last row of
$X$ are zero, and the first row and the first column of $Z$
are zero. Clearly, if $k\ge 1$, then $(X|Z)$ has rank
$n=2k+1$ and $(X,Y)\not\in\overline{G_0(\gt c(1){\times}\gt
c(1))}$. Therefore, for $(\gt{so}_{2n},\gt{gl}_n)$ with odd
$n\ge 3$ condition $({\bf \spadesuit})$ is not satisfied.
\end{ex}

\begin{ex} Consider now symmetric pair
$(E_6,\gt{so}_{10}{\oplus}\bbk)$. Here $\gt g(1)=\bbk^{16}_+$ and 
$\gt g(-1)=\bbk^{16}_-$ are different half-spinor representations of
$\gt{so}_{10}$. Let $\gt t$ be a Cartan subalgebra of $\gt{so}_{10}$,
$\{\pi_i\mid i=1,\ldots, 5\}$  fundamental weights of $\gt{so}_{10}$,
and
$\{\varepsilon_i\mid i=1,\ldots, 5\}$ an orthogonal basis of $\gt t(\mathbb R)^*$
such that $\pi_4=(\varepsilon_1+\varepsilon_2+\varepsilon_3+\varepsilon_4
-\varepsilon_5)/2$, $\pi_5=(\varepsilon_1+\varepsilon_2+\varepsilon_3+\varepsilon_4
+\varepsilon_5)/2$ (for a detailed explanations of this notation see
\cite[Reference Chapter]{vo}).  The rank of $(E_6,\gt{so}_{10}{\oplus}\bbk)$ equals $2$. 
A Cartan subspace $\gt c\subset\gt g_1$ can be chosen such that $\gt c(1)$ is 
a $\gt t$-invariant subspace with weights 
$\pi_5$ and $(\varepsilon_1-\varepsilon_2-\varepsilon_3-\varepsilon_4
-\varepsilon_5)/2$. Let $h\in\gt t$ be an element such that 
$\varepsilon_1(h)=-1$, $\varepsilon_i(h)=0$ for 
$i=2,\ldots, 5$. Clearly, $\lim_{t\in\mathbb Q,\,t\to+\infty}\exp(th){\cdot}v=0$ for 
each $v\in\gt c(1){\times}\gt c(1)\subset\gt g(1){\times}\gt g(1)$.
It remains to find a non-trivial $\SO_{10}$-invariant in 
$\bbk[\gt g(1){\times}\gt g(1)]={\mathcal S}(\bbk^{16}_+{\oplus}\bbk^{16}_-)$.
Denote by $V(\varphi)$ a vector space of a representation with the highest weight 
$\varphi$. Then $\bbk^{16}_+=V(\pi_5)$ and
${\mathcal S}^2(\bbk^{16}_+)=V(2\pi_5)\oplus V(\pi_1)$. 
We have the following  $\SO_{10}$-invariant inclusions:
$$
{\mathcal S^4}(\bbk^{16}_+{\oplus}\bbk^{16}_-)\supset
{\mathcal S^2}(\bbk^{16}_+)\otimes{\mathcal S^2}(\bbk^{16}_-)\supset
 V(\pi_1)\otimes V(\pi_1).
$$
Since $(V(\pi_1){\otimes}V(\pi_1))^{\SO_{10}}=\bbk$, we get a required 
$\SO_{10}$-invariant 
$f$ of degree $4$. Therefore 
$f(\overline{G_0(\gt c(1){\times}\gt c(1)})=0$, 
and condition 
$({\bf \spadesuit})$ is not satisfied. 
\end{ex}

\end{document}